\documentclass[12pt]{amsart}
\usepackage{amsmath}
\usepackage{amsfonts,amssymb,amsthm, txfonts, pxfonts, boxedminipage,
graphicx, color}

\DeclareMathOperator{\med}{med}

\newtheorem{lemma}{Lemma}

\theoremstyle{remark}

\begin{document}

\title{Triangulations into Groups}
\author{Igor Rivin}
\begin{abstract}
If a (cusped) surface $S$ admits an ideal triangulation $\mathcal{T}$ with no
shears, we show an efficient algorithm to give $S$ as a quotient of
hypebolic plane $\mathbb{H}^2$ by a subgroup of $PSL(2, \mathbb{Z}).$
The algorithm runs in time $O(n \log n),$ where $n$ is the number of
triangles in the triangulation $\mathcal{T}.$ The algorithm
generalizes to producing fundamental groups of general surfaces and
geometric manifolds of higher dimension.
\end{abstract}

\date{\today}
\address{Mathematics Department, Temple University, Philadelphia, PA
19122}
\curraddr{Mathematics Department, Fine Hall, Washington Rd, Princeton,
NJ 08544}
\email{irivin@math.princeton.edu}
\thanks{The author would like to thank the Princeton University
Mathematics Department for its hospitality. The first figure was
produced in \textit{Mathematica} using a hyperbolic geometry package
written by the author. The second figure was produced in GraphViz and
massaged by Adobe Illustrator.}

\maketitle
\section{Introduction}

Let $S$ be a cusped hyperbolic surface admitting an ideal
triangulation with the following property: 

\medskip\noindent
\textbf{Property 1.} For any pair of adjacent ideal
triangles $ABC$ and $ABD,$ the \emph{cross ratio} of the four points
$A, B, C, D$ equals $1,$ where the cross ratio is defined as:
\[
\left[A, B, C, D\right] = \dfrac{(A - C)(B - D)}{(B-C)(A-D)},
\]
and we have implicitely identified the hyperbolic plane with the upper
halfplane $\mathbb{H} \subset \mathbb{C}.$

This property has a number of equivalent formulations. One is
geometric:

\medskip\noindent
\textbf{Property 2.}
We can choose a family of horocycles $h_1, \dotsc, h_n,$ where each
$h_i$ is center on the $i$-th cusp of of $S,$ and $h_i$ is tangent to
$h_j$ whenever $c_i$ is adjacent to $c_j$ in $\mathcal{T}.$

Another is algebraic:

\medskip\noindent
\textbf{Property 3.}
$S$ is the quotient of $\mathbb{H}$ by a subgroup $\Gamma$ of $PSL(2,
\mathbb{Z}).$

In this note we will prove that the three properties are equivalent,
and also give an algorithm to express $S = \Gamma \backslash \mathbb{H}.$
The algorithm runs in time bounded by $O(n \log n),$ and produces a
set $G$ of independent matrix generators for $\Gamma.$ Since $S$ is cusped,
$\Gamma$ is a free group, and so this is a complete description of
$\Gamma.$ In fact, we construct the generators as words in the two
linear fractional transformations $L$ and $R,$ where:
\[
L(z) = z + 1, \qquad R(z) = \dfrac{-1}{z-1},
\]
and the matrix generators are obtained by multiplying the words out.

The algorithm consists of a number of steps:

\medskip\noindent
\textbf{Step 1.} Construct the Poincar\'e dual $\mathcal{T}^*$ of the
triangulation. This will have a vertex for each face of $\mathcal{T}$
and a face for each vertex of $\mathcal{T}.$ This is an oriented
complex, and thus we can cyclically order the edges at each vertex.

\medskip\noindent
\textbf{Step 2.} Construct a maximal spanning tree $M$ of the $1$-skeleton
$\mathcal{T}_1^*$ of $\mathcal{T}^*.$ The edges of $\mathcal{T}_1^*$
fall into two types. The edges of the first type. are the edges of
$M,$ the edges of the second type are not. 

\medskip\noindent
\textbf{Step 3.} Split each edge of the second type. By ``split'', I
mean that we replace the edge $AB$ by a pair of edges $AC_1,$
$BC_2.$ We will henceforth refer to $C_1$ and $C_2$ as \emph{twins}.

After we split all the edges of the second type in the graph
$\mathcal{T}_1^*,$ we obtain a graph $B,$ which is a tree where every
non-leaf node has degree three. In addition, every leaf node is
annotated with a cyclic ordering of the three edges. We are ready for:
 
\medskip\noindent
\textbf{Step 4.} Construct the shortest path from each leaf node to
its twin. This path will look like
$C_1 v_1 \dots v_k C_2.$ At each vertex $v_k$ we have a fork in the
road, and we annotate $v_k$ with an $L$ or an $R$ depending on whether
we go left or right at the fork.

Now we are done: each path from $C_1$ to $C_2$ gives a generator of
the fundamental group of $S,$ if we replace $L$ and $R$ by the linear
fractional transformations with the same names (this should be done as
we are constructing the paths, doing it after will bring us back to
$O(n^2)$ running time).

The correctness of the algorithm above follows immediately from the
Poincar\'e Polygon Theorem (see, eg, \cite{beardon}).

\subsection{Crossratios are all $1$ if and only if there is a horodisk
packing}.
This follows from the observation that there is a \emph{unique}
horodisk packing of an ideal triangle. Indeed, if represent the
ideal triangle $ABC$ as one whose vertices are the three roots of unity in the
Poincar\'e disk model, the symmetric arrangement of horodisks
obviously works.Let the points of tangency of the horocycles (which are
on the sides of $ABC$) be $p_{AB},$ $p_{AC},$ and $p_{BC}.$ Now,
suppose that there is another arrangements, with points of tangency
$q_{AB},q_{AC}, q_{BC},$ and let $d(p_{AB}, q_{AB}) = r_{AB},$ and
similarly for the other two sides. Suppose $p_{AB}$ lies between $A$
and $q_{AB}.$ Then the same is true of $p_{AC}$ and $q_{AC}.$ But the
last two assertions would imply the the three new horocycles are not
actually tangent along $BC.$ To show the result we now need the
following easy lemma:

\begin{lemma}
Let $T_1 = ABC$ and $T_2 = ABD$ be two adjacent ideal triangles, and let $\gamma$
be a horocycle centered on $A.$ Let $\gamma_1 = \gamma \cap T_1,$ and
$\gamma_2 = \gamma \cap T_2.$ Then 
\[
\dfrac {|\gamma_1|}{|\gamma_2|} = \exp([A, B, C, D]).
\]
\end{lemma}
\begin{proof}
Let $C = -1,$ $A=\infty,$ $B = 0,$ $D = z,$ and compute.
\end{proof}

\subsection{All crossratios are $1$ implies that the surface is a
quotient of the upper halfplane by a subgroup of the modular group.}
This is not hard.to see, especially if one looks at the modular
figure: each adjacent pair of colored and white triangles forms a
fundamental domain for the action of $PSL(2, \mathbb{Z})$ on the
hyperbolic plane. Since the crossratios are all $1,$ the baricentric
subdivisions of pairs of adjacent ideal triangles agree, and so we
see that our surface covers the modular orbifold.
\begin{figure}
\includegraphics[height=2in]{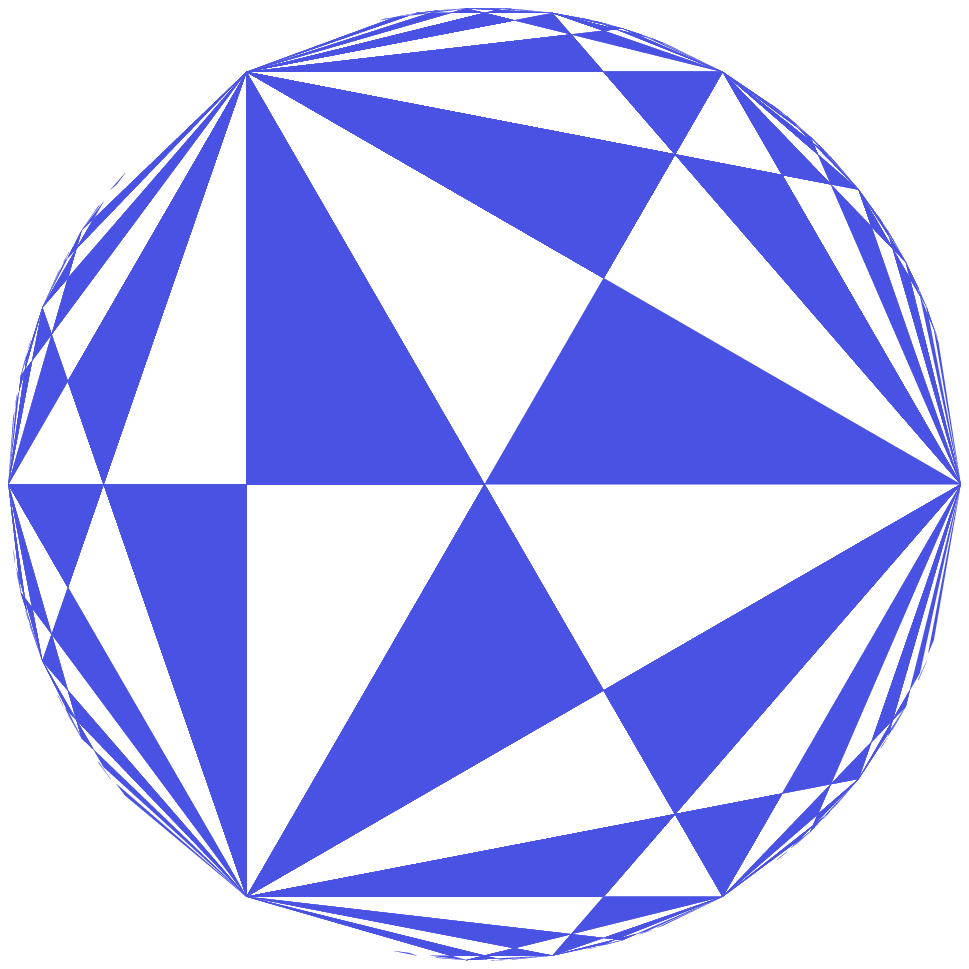}
\caption{The modular tessellation (in the Klein model)}
\end{figure}
\section{Complexity}
\subsection{Constructing the oriented dual (Step 1)}
The complexity of 
Step 1 (constructing the oriented dual) depends on how one is given
the triangulation. The most natural way is for it to be given as a
\emph{rotation system}, which is simply the graph with a cyclic
ordering of the edges at every vertex. It is easy to see that in this
case the dual graph can be constructed in time linear in the number of
edges (the algorithm is simple: maintain a list of edges. Each edge is
marked by 0 or 1. Initially, all the edges have label $0.$ We pick the
first edge $e$, and construct a list of edges obtained by always picking the edge which precedes $e$ in the
cyclic order. A closed cycle gives us a face (already equipped with
the cyclic ordering of boundary edges). Every time an edge is seen we increase the label by $1.$
If the label is $2,$ we delete the edge from the list. Since each
edge is seen at most twice, and we do constant work per edge, the
algorithm is linear). 

The spanning tree (Step 2) can be done in time linear in the number of
edges (see, eg, \cite{cormenetal}), and Step 3 can obviously be done
in time linear in the number of vertices. This leaves us with Step 4,
which we analyze below.

\subsection{Constructing the generators}
At this point we have a tree (with every interior node of degree $3$)
$M$ and a collection of pairs of leaves of $M,$ and we need to
construct paths between the two vertices in each pair. Since a
shortest path between two vertices of a tree can be constructed in
time $O(V(M)),$ (see \cite{cormenetal})  and the number of pairs is half the number of all leaf
nodes (so $O(V(M))$ as well), this gives an $O(V^2(M))$ algorithm for
computing all the generators. We can do better, however, by first
showing the following:
\begin{lemma}
\label{partlem}
Let $G$ be a tree with every non-leaf node having degree three. For
every non-leaf node $v$, removing $v$ separates $G$ into three
subgraphs $G_{\max}(v), G_{\med}(v), G_{\min}(v),.$ with $|V(G_{\max}(v))| \geq
|V(G_{\med}(v))| \geq |V(G_{\min}(v))|.$ Let $\tilde{v}$ be the vertex which
minimizes $|V(G_{\max}(v))|.$ Then
\begin{equation}
\label{lemeq}
\dfrac{V(G)-1}{3}|\leq V(G_{\max}(\tilde{v}))| \leq \frac{2}{3}V(G) +
1.
\end{equation}
\end{lemma}
\begin{proof}
Denote the three orders by $N_{\max},$ $N_{\med},$ $N_{\min}.$
The first inequality is true at any vertex, since 
\begin{equation}
N_{\max} + N_{\med} + N_{\med} + 1 = |V(G)|.
\end{equation}
To show the second inequality, let $v_1$ be the vertex in
$G_{\max}(\tilde{v})$ adjacent to $\tilde{v}.$ At $v_1$ the orders of
the three components into which $v_1$ separates $G$ are $N_1,$ $N_2,$
$N_3,$ where $N_3 = N_{\min} + N_{\med} + 1,$ and $N_1 + N_2 + 1 =
N_{\max}.$ Assume that $N_1 \geq N_2.$ We then have two
possibilities. The first is that $N_3 \geq N_1.$ In that case,
$N_{\med} + N_{\min} + 1 \geq (N_{\max} - 1)/2.$ Adding $N_{\max}$ to
both sides, we get the second side of the inequality.\ref{lemeq}. The
second possibility is that $N_1 \geq N_3.$ However, since $N_1 <
N_{\max},$ this contradicts the defining property of $\tilde{v}.$
\end{proof}
\begin{figure}
\includegraphics[height=2in]{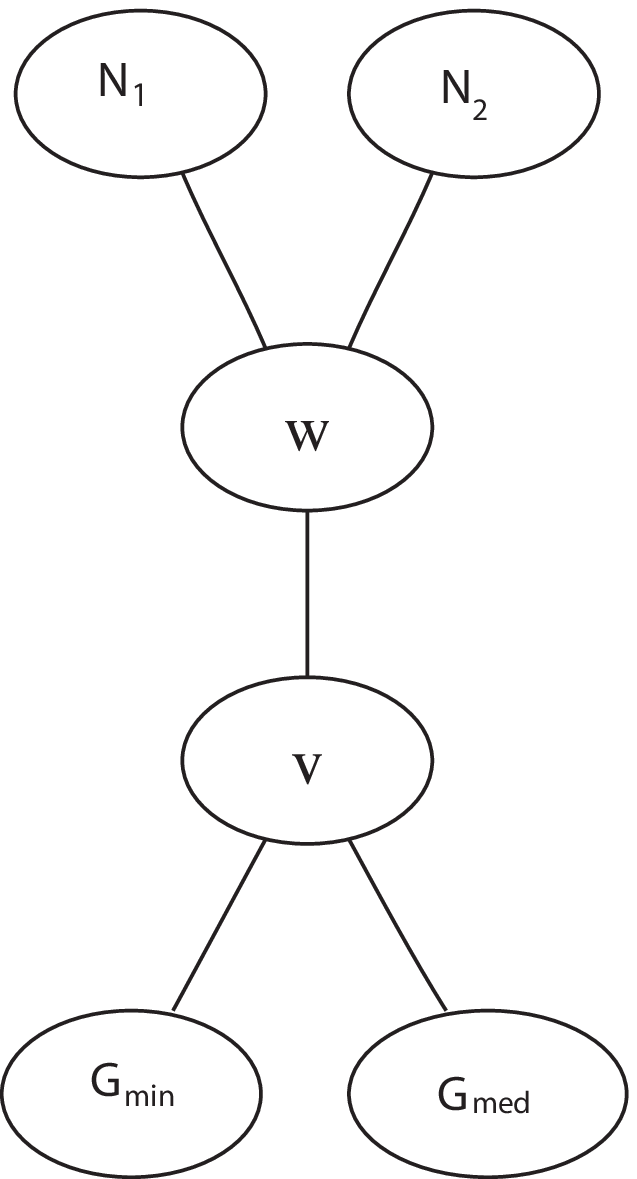}
\caption{Cut point}
\end{figure}

The next result we will need is the result of \cite{kundumisra}:
Given a rooted tree $T$  with a positive weight associated with every node,
there is a linear time algorithm to partition the tree into a minimal
collection of subtrees such that the weight of no subtree exceeds $k.$
In our application, all the weights are equal to $1,$ and $k = \frac23 V(T) + 2.$ By Lemma \ref{partlem} the tree will be broken up
into exactly two components.

The algorithm is then simple: We partiion our tree into two
\emph{rooted} subtrees (the roots will be the two endpoints of the
edge we delete to partition), For each of the two pieces we compute
the pair distances, and all the distances from the leaves to the root
(recursively), then use the distances we had computed to compute all
the distances in the original tree. To make the first step the same as
all the others, we pick an arbitrary leaf and call it the root.

It is clear that at each step we have the following recurrence
inequality for the
running time:

\[
T  \leq c V(G) + T(V_1) + T(V_2),
\]
where $T_1 + T_2 = V(G) - 1,$ and $\max V_1, V_2 \leq \frac23 V(G) + 2.$

This clearly implies an $O(n \log n)$ running time.

\section{Extensions}
The algorithm described above easily extends to other cases. The
simplest extension is where the ideally triangulated surface does not
have all cross ratios equal to one. In that case, we simply replace
the linear fractional transformations $L$ and $R$ by the appropriate
conjugates of the transformation $L_{a, b}(z) = a z + b.$

For a surface (or a higher dimensional manifold) equipped with a triangulation by finite triangles, we
simply develop the fundamental domain (as given by the spanning tree)
into the model space, and then use the side-pairing information to
produce the generators of the fundamental group. It should be noted
that if a surface is finely triangulated, the generating set will
contain many instances of the identity element, but this is obviously
not a serious problem.

It should be noted that Step 4 of our algorithm can be replaced by
using the spanning tree information to embed the triangulation in
$\mathbb{H}^2,$ and then computing the relevant isometries. This
would, however, necessitate running a version of the continued
fraction algorithm for each side pairing, and thus will be much less
efficient ($O(n^2)$ vs our $O(n \log n).$)
\bibliographystyle{plain}
\bibliography{curves,rivin,opt}
\end{document}